# DISCUSSION OF "2004 IMS MEDALLION LECTURE: LOCAL RADEMACHER COMPLEXITIES AND ORACLE INEQUALITIES IN RISK MINIMIZATION" BY V. KOLTCHINSKII

By Gilles Blanchard and Pascal Massart

*Fraunhofer FIRST, Berlin and Université de Paris-Sud*

These last years, much attention has been paid to the construction of model selection criteria via penalization. Vladimir Koltchinskii has to be congratulated for providing a theory reaching a level of generality that is sufficiently high to recover most of the recent results obtained on this topic in the context of statistical learning. Thanks to concentration inequalities and empirical process theory, we are now at a point where the problem of understanding what is the order of the excess risk for the empirical minimizer on a given model is elucidated. Koltchinskii's paper provides several ways of expressing that this excess risk can be sharply bounded by quantities depending on the complexity of the model in various senses. The most prominent relies on Rademacher processes, which Vladimir Koltchinskii himself pioneered in introducing in statistics. We even know that these upper bounds on the excess risk are often unimprovable (see the lower bounds in [6], e.g.).

The same machinery used to analyze the excess risk can be applied to produce penalized criteria and to establish oracle-type risk bounds for the so-defined penalized empirical risk minimizer. The problem of defining properly penalized criteria is particularly challenging in the classification context, since it is connected to the question of defining optimal classifiers *without* knowing in advance the "noise condition" of the underlying distribution [(8.2) of the discussed paper]. This condition determines the attainable rates of convergence and is a topic attracting much attention in the statistical learning community at this moment (see the numerous references in the discussed paper).

What we would like to discuss is the gap between theory and practice of model selection. Of course, the existence of a gap between the methods which are analyzed in theory, and those which are used in practice, is in some sense unavoidable. Our purpose here is to express our perception of the current situation regarding this gap, and to propose some ideas which could contribute to reduce it.

As a starting point for our discussion, we would like to briefly analyze the behavior of the so-called *hold-out* selection procedure. This procedure should be seen as some primitive version of the $V$-fold cross-validation method, which is probably the most commonly used model selection method in practice, in the context of statistical learning. One advantage of hold-out is that it is very easy to study from a mathematical point of view. The point we want to make here is simple, but rich of teachings in this context: hold-out is actually a selection method that is adaptive to the classification noise condition. This property of the hold-out procedure does not seem to be widely known. Since the proof is short and disarmingly simple, we reproduce it here (as inspired from [5], Chapter 8, where more general results for hold-out are also proved).

**1. Hold-out adapts to the noise condition.** Our analysis is based on the following selection theorem among a finite collection of functions, which can be seen as some very elementary and basic version of Theorem 6 of the discussed paper. In what follows, we stick to Vladimir Koltchinskii's notations and conventions; in particular we use (5.3), to express the (unknown) noise conditioning [see also the related equation (8.2)].

1.1. *A basic selection result.*

THEOREM 1. *Let $\{f_m, m \in \mathcal{M}\}$ be a finite collection of real-valued measurable functions defined on some measurable space $\mathcal{X}$ and with $|\mathcal{M}| \geq 2$. Let $\xi_1, \ldots, \xi_n$ be some i.i.d. random variables with common distribution $P$ and denote by $P_n$ the empirical probability measure based on $\xi_1, \ldots, \xi_n$. Assume that $|f_m - f_{m'}| \leq 1$ for every $m, m' \in \mathcal{M}$. Assume furthermore that $Pf_m \geq 0$ for every $m \in \mathcal{M}$.*

*Let $\varphi$ be a convex function on $[0, +\infty)$ with $\varphi(0) = 0$ and such that $\varphi(x)/x^2$ is nondecreasing; denote $\varphi^*$ the convex conjugate of $\varphi$. Assume*

$$\text{(1)} \qquad Pf_m \geq \varphi\left(\sqrt{Pf_m^2}\right) \qquad \text{for every } m \in \mathcal{M}.$$

*Consider some random variable $\widehat{m}$ such that*

$$P_n f_{\widehat{m}} = \inf_{m \in \mathcal{M}} P_n f_m.$$

*Then, for every $\varepsilon \in (0,1)$, the following exponential bound holds for every positive real number $x$:*

$$\text{(2)} \quad \mathbb{P}\left[Pf_{\widehat{m}} > C_\varepsilon \inf_{m \in \mathcal{M}} Pf_m + C'_\varepsilon(x + \ln|\mathcal{M}|)\left(\frac{4}{\varepsilon}\varphi^*\left(\frac{1}{\sqrt{n}}\right) + \frac{1}{3n}\right)\right] \leq e^{-x},$$

where $C_\varepsilon = \frac{1+\varepsilon}{1-\varepsilon}$, $C'_\varepsilon = (1-\varepsilon)^{-1}$.

In particular, the following control in expectation is valid:

$$\text{(3)} \quad \mathbb{E}[Pf_{\widehat{m}}] \leq C_\varepsilon \inf_{m \in \mathcal{M}} \mathbb{E}[Pf_m] + C'_\varepsilon \ln(e|\mathcal{M}|)\left(\frac{4}{\varepsilon}\varphi^*\left(\frac{1}{\sqrt{n}}\right) + \frac{1}{3n}\right).$$

PROOF. Let $m$ be such that
$$Pf_m = \inf_{m' \in \mathcal{M}}(Pf_{m'}).$$

Notice that by definition of $\widehat{m}$, $P_n f_{\widehat{m}} \leq P_n f_m$. Hence,

$$\text{(4)} \quad \begin{aligned} Pf_{\widehat{m}} &= (P - P_n)f_{\widehat{m}} + P_n f_{\widehat{m}} \leq P_n f_m + (P - P_n)(f_{\widehat{m}}) \\ &\leq Pf_m + (P - P_n)(f_{\widehat{m}} - f_m). \end{aligned}$$

Setting for every $m' \in \mathcal{M}$, $\sigma_{m'}^2 = Pf_{m'}^2$, it comes from Bernstein's inequality that for every $m' \in \mathcal{M}$ and every positive number $y$, the following holds except on a set of probability less than $e^{-y}$:

$$(P - P_n)(f_{m'} - f_m) \leq \sqrt{\frac{2y}{n}}(\sigma_m + \sigma_{m'}) + \frac{y}{3n}.$$

By the union bound, choosing $y = \ln|\mathcal{M}| + x$, and using (1), this implies that, except on some set $\Omega_x$ with probability less than $e^{-x}$,

$$\text{(5)} \quad (P - P_n)(f_{\widehat{m}} - f_m) \leq \sqrt{\frac{2y}{n}}(\varphi^{-1}(Pf_{\widehat{m}}) + \varphi^{-1}(Pf_m)) + \frac{(\ln|\mathcal{M}| + x)}{3n}.$$

Let $\varphi^*$ be the convex conjugate of $\varphi$; we then have

$$\sqrt{\frac{2y}{n}}(\varphi^{-1}(Pf_m)) \leq \varphi(\sqrt{\varepsilon}\varphi^{-1}(Pf_m)) + \varphi^*\left(\sqrt{\frac{2y}{\varepsilon n}}\right) \leq \varepsilon Pf_m + \frac{2y}{\varepsilon}\varphi^*\left(\frac{1}{\sqrt{n}}\right),$$

with a similar inequality for $\widehat{m}$. For the last inequality above, we have used the assumption that $\varphi(x)/x^2$ is nondecreasing, which readily implies that $\varphi^*(x)/x^2$ is nonincreasing, along with the fact that $\varepsilon \leq 1$ and $2y/\varepsilon \geq 1$. Combining this inequality with (5) and (4) yields

$$(1-\varepsilon)P(f_{\widehat{m}}) \leq (1+\varepsilon)P(f_m) + (x + \ln|\mathcal{M}|)\left(4\varepsilon^{-1}\varphi^*\left(\frac{1}{\sqrt{n}}\right) + \frac{1}{3n}\right). \quad \square$$

1.2. *An oracle inequality for hold-out.* Let us now describe and study the hold-out procedure. Assume that we observe $N + n$ random variables with common distribution $P$ depending on some parameter $g_*$ to be estimated. The first $N$ observations $\xi'_1, \ldots, \xi'_N$ are used to build some preliminary collection of estimators $\{\widehat{g}_m\}_{m \in \mathcal{M}}$ and we use the remaining observations $\xi_1, \ldots, \xi_n$ to select some estimator $\widehat{g}_m$ among the collection $\{\widehat{g}_m\}_{m \in \mathcal{M}}$. We more precisely consider here the situation described in Section 7 of the paper, where there is some (bounded) loss or contrast

$$\ell : T \times \mathbb{R} \to [0,1]$$

which is well adapted to our estimation problem of $g_*$ in the sense that the expected loss $\mathbb{E}[\ell \bullet g] = \mathbb{E}[\ell(Y, g(X))]$ achieves a minimum at $g_*$ when $g$ varies in $\mathcal{G}$. We denote the relative expected loss as follows:

$$L(g, g_*) = \mathbb{E}[\ell \bullet g - \ell \bullet g_*] \qquad \text{for all } g \in \mathcal{G}.$$

For bounded regression or binary classification, we can take, for example,

$$\ell(y, x) = (y - x)^2;$$

then $g_*(x) = P(Y = 1 | X = x)$ (resp. $g_* = b^*$, the Bayes classifier) is indeed the minimizer of $\mathbb{E}[(Y - t(X))^2]$ over the set of measurable functions $g$ taking their values in $[0,1]$ (resp. $\{0,1\}$). We can now apply Theorem 1, conditionally on the training sample $\xi_1', \ldots, \xi_N'$, to the collection of functions

$$\{f_m = (\ell \bullet \widehat{g}_m - \ell \bullet g_*), m \in \mathcal{M}\}.$$

Define, as in the theorem, $\widehat{m}$ as a minimizer of the empirical risk $P_n(\ell \bullet \widehat{g}_m)$ over $\mathcal{M}$. If $\varphi$ satisfies the weak regularity assumptions of Theorem 1 and is such that

(6) $$\sup_{L(g, g_*) \leq \varepsilon} \|\ell \bullet g - \ell \bullet g_*\|_{2, P} \leq \varphi^{-1}(\varepsilon),$$

we derive from (3) that conditionally on $\xi_1', \ldots, \xi_N'$, one has for every $\varepsilon \in (0, 1)$:

(7) $$\mathbb{E}[L(\widehat{g}_{\widehat{m}}, g_*) | \xi'] \leq C_\epsilon \inf_{m \in \mathcal{M}} L(\widehat{g}_m, g_*) + C'_\varepsilon \ln(e|\mathcal{M}|) \left( \frac{4}{\varepsilon} \varphi^* \left( \frac{1}{\sqrt{n}} \right) + \frac{1}{3n} \right).$$

The striking feature of this result is that the hold-out selection procedure provides an oracle-type inequality involving the modulus of continuity $\varphi^{-1}$ which is not known in advance. This is especially interesting in the classification framework for which $\varphi$ can be of very different natures according to the difficulty of the classification problem. The main issue is therefore to understand whether the term $\varphi^*(n^{-1/2})(1 + \ln|\mathcal{M}|)$ appearing in (7) is indeed a remainder term or not. We cannot exactly answer this question in general because it is hard to compare $\delta_n := \varphi^*(n^{-1/2})$ with $\inf_{m \in \mathcal{M}} L(\widehat{g}_m, g_*)$. However, if $\widehat{g}_m$ is itself an empirical risk minimizer over some model $\mathcal{G}_m$, we can compare $\delta_n$ with $\inf_{m \in \mathcal{M}} \theta_{m,N}$, where $\theta_{m,N}$ is an upper bound (up to constant) for the expectation of the excess risk within model $\mathcal{G}_m$. More precisely, taking for instance $\varepsilon = 1/2$, we derive from (7) that for some constant $\kappa$

(8) $$\mathbb{E}[\ell(\widehat{g}_{\widehat{m}}, g_*)] \leq 3 \inf_{m \in \mathcal{M}} (L(\mathcal{G}_m, g_*) + \kappa \theta_{m,N}) + \ln(e|\mathcal{M}|)(16\delta_n + (3n)^{-1}),$$

where $L(\mathcal{G}_m, g_*) = \inf_{g \in \mathcal{G}_m} L(g, g_*)$. Now, using the notation and method described in Section 3 of the discussed paper, we obtain as a value for $\theta_{m,N}$ the (largest) fixed point of the function

$$\bar{U}_N^{(m)}(\delta) = K(\phi_{m,N}(\delta) + D(\delta)N^{-1/2} + N^{-1}),$$

where $\phi_{m,N}$ is a nondecreasing function which more or less plays the role of a modulus of continuity of the empirical process $(P'_N - P)(\ell \bullet g)$ over the model $\mathcal{G}_m$. If $N$ and $n$ are of the same order of magnitude, say $N = n$ to be as simple as possible, then $\bar{U}_n^{(m)}(\delta) \geq KD(\delta)n^{-1/2} \geq K'\varphi^{-1}(\delta)n^{-1/2}$, where we assumed that (6) is sharp up to a constant. Therefore, $\theta_{m,n}$ is surely larger (again up to constant factor) than the solution $\widetilde{\delta}_n$ of $\varphi(\delta) = \delta n^{-1/2}$. Since $\varphi$ is a nonnegative convex function with $\varphi(0) = 0$, elementary considerations then show that $\widetilde{\delta}_n \geq \delta_n$.

In fact, $\theta_{m,n}$ will typically turn up much larger in magnitude (as a function of $n$) than $\delta_n$, since the fixed point equation for $\theta_{m,n}$ also involves the function $\phi_{m,n}$, which measures in some way the complexity of model $m$. Finally, the factor $\ln |\mathcal{M}|$ appearing above should also be considered minor if we assume, for instance, that $|\mathcal{M}| \leq n^k$ for some $k \geq 0$. In conclusion, except in very pathological situations, the quantity $\delta_n \ln |\mathcal{M}|$ really plays the role of a remainder term in (8).

## 2. Data-driven penalties.

2.1. *A sober assessment of the current state of the art.* It is, in some sense, somewhat disappointing to discover that a very crude method like hold-out is working so well. This is especially true in the classification framework, where it is indeed painstakingly difficult to design penalties that are adaptive to the noise condition. Recent works on the topic, involving local Rademacher penalties for instance, provide at least some theoretical solutions to the problem; but they systematically involve unknown constants—either because the numerical values coming from the theory are overpessimistic, or, worse, because these constants also depend on nuisance parameters related to the unknown distribution (e.g., the infimum of the density of explanatory variables).

We therefore end up in the following delicate situation:

- From a *theoretical* point of view, we are not in a position to justify that conveniently penalized model selection methods (or more generally, model selection methods that use the entire sample for the estimation within each model) could improve over the simple hold-out solution.
- From a *practical* point of view, the penalization method does not provide a "ready-to-use" solution and remains far from being competitive with relatively simple methods that are widely used in practice. We have in mind in particular $V$-fold cross-validation.

At this stage, two natural and connected questions emerge:

- Is there some room left for penalization methods?
- How should penalties be calibrated to design efficient procedures?

There is at least one strong reason for which, despite the arguments developed above, one should keep interest in penalization methods: for independent but not identically distributed observations (we typically think of Gaussian regression on a fixed design), hold-out (for theory) or cross-validation (for practice) may break or become irrelevant.

Another issue is that, intuitively, one would expect that the hold-out leads to the loss of a factor 2 because of the sample halving process. Unfortunately, much looser constants appear when applying the more elaborate theoretical tools needed to tackle penalization, where the entire sample is used for the estimation within each model. These larger constants drown out the initial factor 2 advantage over the hold-out, so that the theory may currently not be precise enough to distinguish this effect.

In other words, since the opponents are strong, beating them remains possible but requires one to calibrate penalties sharply. This leads us to the second question raised above. We would like now to provide some ideas that can contribute to answering this last question, partly based on theoretical results which are already available and partly based on heuristics and thoughts which lead to some empirical rules and new theoretical problems.

2.2. *A rule of thumb for calibrating penalties from the data.* A general idea consists in guessing what is the right penalty to be used from the data itself. Let us roughly describe the type of results which been proved in the Gaussian framework in [2]. In several contexts (such as variable selection, e.g.), it is possible to prove lower bounds for penalties (meaning that lower penalties will lead to asymptotic inconsistency). Moreover, a close inspection of oracle inequalities shows that approximately optimal values for the penalty are linked to minimal values within a factor 2. We can therefore retain from this Gaussian theory the rule of thumb:

$$\text{(9)} \qquad \text{``optimal'' penalty} = 2 \times \text{``minimal'' penalty}.$$

Interestingly, the *minimal* penalty can be evaluated from the data: when the penalty is not heavy enough, one systematically chooses models with very large dimension. It remains to double this minimal penalty to produce the desired (nearly) optimal one. This strategy allows to design a data-driven penalty without knowing in advance the level of noise in Gaussian regression. In the context of change points detection, this data-driven calibration method for the penalty has been successfully implemented and tested by Lebarbier (see [3]).

In the non-Gaussian case, we believe that this procedure remains valid, but theoretical justification remains an open problem. As already mentioned earlier, this problem is especially challenging in the classification context, since it is connected to the question of defining optimal classifiers *without* knowing in advance the noise condition of the underlying distribution.

2.3. *Akaike's heuristics revisited.* In order to better understand the above rule of thumb and understand why it could be extended to non-Gaussian frameworks, it is instructive to come back to the original ideas of model selection via penalization, that is, Mallows's or Akaike's heuristics (see [4] and [1]). Both are based on the principle of unbiased estimation of the risk (at least asymptotically as far as Akaike's heuristics is concerned). Our idea is to adapt this principle to a nonasymptotic view of the question, for which one could hope to use concentration inequalities rather than limit theorems to validate the heuristics.

Let us consider, in each model $\mathcal{G}_m$, some minimizer $g_m$ of $g \to \mathbb{E}[\ell \bullet g]$ over $\mathcal{G}_m$ (assuming that such a point does exist). Defining for every $m \in \mathcal{M}$,

$$\widehat{b}_m = P_n(\ell \bullet g_m - \ell \bullet g_*) \quad \text{and} \quad \widehat{v}_m = P_n(\ell \bullet g_m - \ell \bullet \widehat{g}_m),$$

minimizing some penalized criterion $P_n(\ell \bullet \widehat{g}_m) + \text{pen}(m)$ over $\mathcal{M}$ amounts to minimizing $\widehat{b}_m - \widehat{v}_m + \text{pen}(m)$. The point is that $\widehat{b}_m$ is an unbiased estimator of the bias term $L(g_m, g_*)$. With concentration arguments in mind, one can hope that minimizing the above quantity will be approximately equivalent to minimizing $L(g_m, g_*) - \mathbb{E}[\widehat{v}_m] + \text{pen}(m)$. Since the purpose of the game is to minimize the risk $\mathbb{E}[L(\widehat{g}_m, g_*)]$, an ideal penalty would therefore be

$$\text{pen}(m) = \mathbb{E}[\widehat{v}_m] + \mathbb{E}[L(\widehat{g}_m, g_m)].$$

In Mallows's $C_p$ case, $\ell$ is the square loss, the models $\mathcal{G}_m$ are linear and $\mathbb{E}[\widehat{v}_m] = \mathbb{E}[L(\widehat{g}_m, g_m)]$ are explicitly computable (at least if the level of noise is assumed to be known). For Akaike's penalized log-likelihood criterion, this is similar, at least asymptotically. More precisely, in Akaike's heuristics, $\ell$ is the (minus) log-likelihood and one uses the fact that $\mathbb{E}[\widehat{v}_m] \approx \mathbb{E}[L(\widehat{g}_m, g_m)] \approx D_m/(2n)$, where $D_m$ stands for the number of parameters defining model $\mathcal{G}_m$.

Of course, we do not want to take in consideration the second approximation, which is typically asymptotic and relies on the specific choice of the log-likelihood loss, as well as on regularity conditions of the parametric models, that we certainly do not want to assume here. Our guess, however, is that one can trust the first approximation $\mathbb{E}[\widehat{v}_m] \approx \mathbb{E}[L(\widehat{g}_m, g_m)]$ in a more general situation. If one believes in the validity of this approximation, then a good penalty is $2\mathbb{E}[\widehat{v}_m]$, or equivalently (having still in mind concentration arguments) $2\widehat{v}_m$. This, in some sense, explains the rule of thumb which is given in the preceding section: the minimal penalty is $\widehat{v}_m$, while the optimal penalty should be $\widehat{v}_m + \mathbb{E}[L(\widehat{g}_m, g_m)]$, and their ratio is approximately equal to 2 [note that Akaike's criterion itself may be interpreted by formula (9), the minimal penalty being taken as $D_m/(2n)$ and the optimal penalty as $D_m/n$]. As mentioned above, the interesting point is that, even though $\widehat{v}_m$ is not observable, we can guess the minimal penalty from the data anyway. One way to do this in practice is to search a minimal penalty of the form $\text{pen}(m) = \alpha D_m$ and estimate $\alpha$ by choosing the smallest value for which the corresponding penalized criterion does not lead to selecting "very large" models. Of course, concentration arguments will work only if the list of models is not too rich. In practice, this means that, starting from a given list of models, one has first to decide to penalize in the same way the models which are defined by the same number of parameters. Then one considers a new list of models $(\mathcal{G}_D)_{D \geq 1}$, where for each integer $D$, $\mathcal{G}_D$ is the union of those among the initial models which are defined by $D$ parameters and then applies the preceding heuristics to this new list.

**3. Conclusion.** Caricaturing a little, we could say that at this point, we have a beautiful, yet not very useful, theory—at least, this is the conclusion to which a person mainly interested in practical applications could come. Hold-out is our nemesis for theory, as is cross-validation for practice. Moreover, note that hold-out is also known to be quite unstable in practice—this is the reason why cross-validation is preferred—which widens the gap theory/practice yet a little more.

An optimistic way to look at this, though, is to say that we are only half-way climbing the slope, and that many interesting problems are open to future research efforts. Obviously, there are at least two directions of research. The first one consists in designing proper data-driven penalties which are ready to be used in practice and theoretically efficient (we have tried to give some ideas in this direction in the preceding section). In the spirit of the above results on hold-out, the second one consists in studying in depth the theoretical properties of $V$-fold cross-validation.

Fraunhofer FIRST.IDA
Kekuléstrasse 7, 12489 Berlin
Germany
E-mail: blanchar@first.fhg.de

Laboratoire de mathématiques
Université Paris-Sud – Bât. 425
91405 Orsay Cedex
France
E-mail: Pascal.Massart@math.u-psud.fr